\documentclass[12pt]{article}
\usepackage{amsfonts,amssymb}
\usepackage{amsmath}
\usepackage{amsthm}
\newtheorem*{theorem}{Theorem}
\newtheorem{lemma}{Lemma}

\newtheorem{corollary}{Corollary}
\newtheorem{definition}{Definition}
\title{Families of Vector Fields which Generate the Group of Diffeomorphisms}
\author{A.~A.~Agrachev\thanks{SISSA, Trieste \& Steklov Math. Inst.,
Moscow} \and M.~Caponigro\thanks{SISSA, Trieste}}
\date{}
\begin{document}
\maketitle
\begin{abstract} Given a compact manifold $M$, we prove that any
bracket generating and invariant under multiplication on smooth
functions family of vector fields on $M$ generates the connected
component of unit of the group $\mathrm{Diff}M$.
\end{abstract}

Let $M$ be a smooth\footnote{In this paper, smooth means
$C^\infty$.} $n$-dimensional compact manifold, $\mathrm{Vec}M$ the
space of smooth vector fields on $M$ and $\mathrm{Diff}_0M$ the
group of isotopic to the identity diffeomorphisms of $M$.

Given $f\in\mathrm{Vec}M$, we denote by $t\mapsto e^{tf},\
t\in\mathbb R,$ the flow on $M$ generated by $f$; then $e^{tf},\
t\in\mathbb R,$ is a one-parametric subgroup of
$\mathrm{Diff}_0M$. Let $\mathcal F\subset\mathrm{Vec}M$; the
subgroup of $\mathrm{Diff}_0M$ generated by $e^{tf},\ f\in\mathcal
F,\ t\in\mathbb R,$ is denoted by $\mathrm{Gr}\mathcal F$.

\begin{theorem} Let $\mathcal F\subset \mathrm{Vec}M$; if
$\mathrm{Gr}\mathcal F$ acts transitively on $M$, then
$$
\mathrm{Gr}\left\{af:a\in C^\infty(M),\ f\in\mathcal
F\right\}=\mathrm{Diff}_0M.
$$
\end{theorem}

\begin{corollary} Let $\Delta\subset TM$ be a completely nonholonomic
vector distribution. Then any isotopic to the identity
diffeomorphism of $M$ has a form $e^{f_1}\circ\cdots\circ
e^{f_k}$, where $f_1,\ldots,f_k$ are sections of $\Delta$.
\end{corollary}

\noindent{\bf Remark.} Recall that $\mathrm{Gr}\{f_1,f_2\}$ acts
transitively on $M$ for a generic pair of smooth vector fields
$f_1,f_2$.

\bigskip We start the proof of the theorem with an auxiliary lemma
that is actually the main part of the proof. Let $B\subset\mathbb R^n$
be diffeomorphic to a cube, $0\in B$; we set
$C_0^\infty(B)=\{a\in C^\infty(B):a(0)=0\}$
and assume that $C^\infty_0(B)$ is endowed with the
standard $C^\infty$-topology.
\begin{lemma}[Main Lemma] Let $X_i\in\mathrm{Vec}\mathbb R^n,\
a_i\in C^\infty(\mathbb R^n),\ i=1,\ldots,n$, and the
following conditions hold:
\begin{itemize} \item
$span\{X_1(0),\ldots,X_n(0)\}=\mathbb R^n$,
\item $a_i(0)=0,\ \langle d_0a_i,X_i(0)\rangle<0,\ i=1,\ldots,n$;
\end{itemize}
then there exist $\epsilon,\varepsilon>0$ and a neighborhood
$\mathcal O$ of $(\epsilon a_1,\ldots,\epsilon
a_n)\bigr|_{B_\varepsilon}$ in $C_0^\infty(B_\varepsilon)^n$ such
that the mapping
$$
\Phi:(b_1,\ldots,b_n)\mapsto \left(e^{b_1X_1}\circ\cdots\circ
e^{b_nX_n}\right)\bigr|_{B_{\varepsilon}} \eqno (1)
$$
is an open map from $\mathcal O$ into
$C_0^\infty(B_{\varepsilon})^n$, where
$$
B_\varepsilon=\left\{e^{s_1X_1}\circ\cdots\circ e^{s_nX_n}(0):
|s_i|\le\varepsilon,\ i=1,\ldots,n\right\}.
$$
\end{lemma}
{\bf Sketch of proof.} Openness of the map (1) is derived from the
Hamilton's version of the Nash--Moser inverse function theorem
\cite{Ha}. Set $\bar a=(\epsilon a_1,\ldots,\epsilon a_n)$. In
order to apply the Nash--Moser theorem we have to invert the
differential of $\Phi$ at $\bar a$ and show that inverse is
``tame" with respect to $\bar a$. Here we make computations only
for fixed $\bar a$ and leave the boring check of the tame
dependence on $\bar a$ for the detailed paper.

Note that $e^{\epsilon a_jX_j}$ are closed to identity
diffeomorphisms, hence $\frac{\partial\Phi}{\partial
b_i}\bigr|_{\bar a}$ is obtained from $\frac{\partial}{\partial
b_i}e^{b_iX_i}\bigr|_{\epsilon a_i}$ by a closed to identity
change of variables. We have
$$
\left(\frac{\partial}{\partial a}e^{aX}\right):b\mapsto e^{aX}_*
\left(\int_0^1 e^{\int_t^0\langle da,X\rangle\circ e^{\tau
aX}d\tau}b\circ e^{taX} \,dtX\right)\circ e^{aX}.
$$
This equality follows from the standard ``variations formula" (see
\cite{AS}) and the relation:
$$
\left(e^{taX}\right)_*:X\mapsto \left(e^{\int_0^t\langle
da,X\rangle\circ e^{-\tau aX}d\tau}\right)X.
$$

Let us define an operator
$A(a,X):C_0^\infty(\hat B_\varepsilon)\to C_0^\infty\left(\hat B_{\varepsilon}\right)$
by the formula
$$
A(a,X)b=\int_0^1 e^{\int_t^0\langle da,X\rangle\circ e^{\tau
aX}d\tau}b\circ e^{taX}\,dt,
$$
where $\hat B_{\varepsilon}=\left\{e^{sX}(x):|s|\le\varepsilon,\
x\in\Pi^{n-1}\right\}$ and $\Pi^{n-1}$ is a transversal to $X$
small $(n-1)$-dimensional box. We see that invertibility of
$A(\varepsilon a_i,X_i),\ i=1,\ldots,n,$ implies invertibility of
$D_{\bar a}\Phi$.

Now set $\mathcal X=\{bX: b\in C^\infty(M)\}\subset\mathrm{Vec}M$.
The map
$$
(bX)\mapsto\left(A(a,x)b\right)X
$$
has a clear intrinsic meaning as a linear operator on the space
$\mathcal X$; moreover, this operator depends only on the vector
field $aX\in\mathcal X$. Indeed,
$$
\left(A(a,X)b\right)X=e^{-aX}_*\left(D_{(aX)}Exp\bigr|_{\mathcal
X}(bX)\right)\circ e^{-aX},
$$
where $D_YExp$ is the differential at the point
$Y\in\mathrm{Vec}M$ of the map
$$
Exp:Y\mapsto e^Y,\quad Y\in\mathrm{Vec}M.
$$

Recall that $a(0)=0,\ \langle d_0a,X(0)\rangle<0$. In particular,
$X$ is transversal to the hypersurface $a^{-1}(0)$. We may rectify
the field $X$ in such a way that, in new coordinates,
$X=\frac\partial{\partial x_1},\ a(0,x_2,\ldots,x_n)=0$. Now the
field $aX$ can be treated as a depending on $y=(x_2,\ldots,x_n)$
family of 1-dimensional vector fields
$a(x_1,y)\frac\partial{\partial x_1}$. Moreover,
$a(0,y)=0,\quad \frac{\partial a}{\partial x_1}(0,y)=\alpha(y)<0.$

A hyperbolic 1-dimensional field $a(x_1,y)\frac\partial{\partial
x_1}$ can be linearized by a smooth change of variable and this
smooth change of variable smoothly depends on $y$. Hence we may
assume that $aX=\alpha(y)x_1\frac\partial{\partial x_1}$. Then
$b\circ e^{taX}(x_1,y)=b(e^{\alpha(y)t}x_1,y)$.

We thus have to invert the operator
$$
\hat A:
b(x_1,y)\mapsto\int\limits_0^1e^{-t\alpha(y)}b\left(e^{\alpha(y)t}x_1,y\right)\,dt
$$
acting in the space of smooth functions on a box. We can write
$$b(x_1,y)=b_0(y)+x_1b_1(y)+x^2_1u(x_1,y),$$ where $u$ is a smooth function.
Then $\hat Ab_0=\frac 1\alpha(1-e^{-\alpha})b_0,\ \hat
A\left(x_1b_1\right)=x_1b_1$ and
$$
\hat A\left(x_1^2u(x_1,y)\right)=
x^2_1\int\limits_0^1e^{\alpha(y)t}u\left(e^{\alpha(y)t}x_1,y\right)\,dt
= -\frac{x_1^2}{\alpha(y)}\int\limits_{e^{\alpha(y)}}^1u(\tau
x_1,y)\,d\tau.
$$

What remains is to invert the operator
$$
B:u(x_1,y)\mapsto \int_{e^{\alpha(y)}}^1u(\tau x_1,y)\,d\tau.
$$
We set $v(x_1,y)=\frac 1{x_1}\int_0^{x_1}u(s,y)\,ds$; then
$$
(Bu)(x_1,y)=\left(v(x_1,y)-e^{\alpha(y)}v(e^{\alpha(y)}x_1,y)\right).
\eqno (2)
$$
We introduce one more operator:
$$
R:v(x_1,y)\mapsto e^{\alpha(y)}v(e^{\alpha(y)}x_1,y).
$$
Let $\|v\|_{C^{k,0}}= \sup\limits_{1\le i\le
k}\bigr\|\frac{\partial^iv}{\partial x_1^i}\bigr\|_{C^0}.$
Obviously, $\|R\|_{C^{k,0}}\le e^{\sup\alpha}<1,\ \forall k$.
Hence $(I-R)^{-1}$ transforms a smooth on the box function $\psi$
in the function $\varphi=(I-R)^{-1}\psi$ that is smooth with
respect to $x_1$. As usually, the chain rule for the
differentiation allows to demonstrate that function $\varphi$ is
also smooth on the box and to compute its derivatives:
$$
\frac{\partial\varphi}{\partial y_i}=\left(I-R\right)^{-1}\left(
\frac{\partial\psi}{\partial y_i}-e^\alpha\frac{\partial\alpha}{\partial y_i}\varphi
-e^{2\alpha}\frac{\partial\alpha}{\partial y_i}\frac{\partial\varphi}{\partial x_1}\right),
\quad \mathrm{e. t. c.}
$$
Coming back to equation (2), we obtain: $v=(I-R)^{-1}Bu$. Finally,
$$
B^{-1}:w\mapsto\frac\partial{\partial x_1}\left(x_1(I-R)^{-1}
w\right). \eqno \square
$$

Now set
$$
\mathcal P=\mathrm{Gr}\left\{af:a\in C^\infty(M),\ f\in\mathcal
F\right\}, \quad \mathcal P_q=\{P\in\mathcal P:P(q)=q\},\ q\in M.
$$

\begin{lemma} Any $q\in M$ possesses a neighborhood $U_q\subset M$
such that the set
$$
\left\{P\bigr|_{U_q}:P\in\mathcal P_q\right\} \eqno (3)
$$
has a nonempty interior in $C_q^\infty(U_q,M)$, where
$C_q^\infty(U_q,M)$ is the Fr\'echet manifold of smooth maps
$F:U_q\to M$ such that $F(q)=q$.
\end{lemma}
{\bf Proof.} According to the Orbit Theorem of Sussmann \cite{Su}
(see also the textbook \cite{AS}),
transitivity of the action of $\mathrm{Gr}\mathcal F$ on $M$
implies that
$$
T_qM=span\{P_*f(q):p\in\mathrm{Gr}\mathcal F,\ f\in\mathcal F\}.
$$
Take $X_i={P_i}_*f_i,\ i=1,\ldots,n,$ such that
$P_i\in\mathrm{Gr}\mathcal F,\ f_i\in\mathcal F$, and
$X_1(q),\ldots,X_n(q)$ form a basis of $T_qM$. Then for any
vanished at $q$ smooth functions $a_1,\ldots,a_n$, the
diffeomorphism
$$
e^{a_1X_1}\circ\cdots\circ e^{a_nX_n}=P_1\circ e^{(a_1\circ
P_1)f_1}\circ P^{-1}_1\circ\cdots\circ P_n\circ e^{(a_n\circ
P_n)f_n}\circ P^{-1}_n
$$
belongs to the group $\mathcal P_q$. The desired result now
follows from Main Lemma.

\begin{corollary} Interior of the set (3) contains the identical
map.
\end{corollary}
{\bf Proof.} Let $\mathcal O$ be an open subset of
$C_q^\infty(U_q,M)$ that is contained in (3) and
$P_0\bigr|_{U_q}\in\mathcal O$. Then $P_0^{-1}\circ\mathcal O$ is
a contained in (3) neighborhood of the identity.

\begin{definition} Given $P\in\mathrm{Diff}M$, we set
$supp\,P=\overline{\{x\in M:P(x)\ne x\}}$.
\end{definition}

\begin{lemma} Let $\mathcal O$ be a neighborhood of the identity
in $\mathrm{Diff}M$. Then for any $q\in M$ and any neighborhood
$U_q\subset M$ of $q$, we have:
$$
q\in int\left\{P(q):P\in\mathcal O\cap\mathcal P,\ supp\,P\subset
U_q\right\}.
$$
\end{lemma}
{\bf Proof.} Let vector fields $X_1,\ldots, X_n$ be as in the
proof of Lemma~2 and $b\in C^\infty(M)$ a cut-off function such
that $supp\,b\subset U_q$ and $q\in int\, b^{-1}(1)$. Then the
diffeomorphism
$$
Q(s_1,\ldots,s_n)=e^{s_1bX_1}\circ\cdots\circ e^{s_nbX_n}
$$
belongs to $\mathcal O\cap\mathcal P$ for all sufficiently close
to 0 real numbers $s_1,\ldots,s_n$ and
$supp\,Q(s_1,\ldots,s_n)\subset U_q$. On the other hand, the map
$$(s_1,\ldots,s_n)\mapsto Q(s_1,\ldots,s_n)(q)$$ is a local
diffeomorphism in a neighborhood of 0.

\begin{lemma} Let $\bigcup\limits_jU_j=M$ be a covering of $M$ by
open subsets and $\mathcal O$ be a neighborhood of identity in
$\mathrm{Diff}M$. Then the group $\mathrm{Diff}_0M$ is generated
by the subset
$$
\{P\in\mathcal O: \exists j\ \mathrm{such\ that}\ supp\,P\subset
U_j\}.
$$
\end{lemma}
{\bf Proof.} The group $\mathrm{Diff}_0M$ is obviously generated
by any neighborhood of the identity. We may assume that the
covering of $M$ is finite and any $U_j$ is contained in a
coordinate neighborhood. Moreover, taking a finer covering and a
smaller neighborhood $\mathcal O$ if necessary, we may assume that
for any $P\in\mathcal O$ and any $U_j$, the coordinate
representation of $P\bigr|_{U_j}$ has a form $P:x\mapsto
x+\varphi_P(x)$, where $\varphi$ is a $C^1$-small smooth vector
function.

Now consider a refined covering $\bigcup\limits_iO_i=M$, so that
$\overline{O}_i\subset U_{j_i}$ for some $j_i$ and cut-off
functions $a_i$ such that $a_i|_{O_i}=1,\ supp\,a_i\subset
U_{j_i}$. Given $P\in\mathcal O$, we set
$$
P_i(x)=x+a_i(x)\varphi_P(x),\ \forall x\in U_{j_i}\ \mathrm{and}\
P_i(q)=q,\ \forall q\in M\setminus U_{j_i}.
$$
Then $supp\,(P_i^{-1}\circ P)\subset supp\,P\setminus O_i$. Now,
by the induction with respect to $i$, we step by step arrive to a
diffeomorphism with empty support. In other words, we present $P$
as a composition of diffeomorphisms whose supports are contained
in $U_j$.

\medskip\noindent{\bf Proof of the Theorem.} According to Lemma~4,
it is sufficient to prove that there exist a neighborhood
$U_q\subset M$ and a neighborhood of the identity $\mathcal
O\subset\mathrm{Diff}M$ such that any diffeomorphism $P\in\mathcal
O$ whose support is contained in $U_q$ belongs to $\mathcal P$.
Moreover, Lemma~3 allows to assume that $P(q)=q$. Finally, the
corollary to Lemma~2 completes the job.

\medskip\noindent{\sl Acknowledgment.} First coauthor is greatful to
Boris Khesin who asked him the question answered by this paper
(see also recent preprint \cite{KhL}).


\begin{thebibliography}{9}

\bibitem{AS} A.~A.~Agrachev,Yu.~L.~Sachkov, {\it Control theory from
the geometric viewpoint}. Springer-Verlag, Berlin, 2004,
xiv+412pp.

\bibitem{Ha} R. Hamilton, {\it The inverse function theorem of Nash
and Moser}. Bulletin of the Amer. Math. Soc., 1982, v.7, 65--222

\bibitem{KhL} B. Khesin, P. Lee, {\it A nonholonomic Moser theorem
and optimal mass transport}. ArXiv:0802.1551v2 [math.DG], 2008,
31p.

\bibitem{Su} H. J. Sussmann, {\it Orbits of families of vector fields
and integrability of distributions}. Trans. Amer. Math. Soc., 1973,
v.180, 171--188
\end{thebibliography}
\end{document}